\documentclass[10pt]{amsart}

\usepackage{amsfonts,graphics,amsmath,amsthm,amsfonts,amscd,
amssymb,amsmath,latexsym}

\usepackage{epsfig}

\theoremstyle{plain}

\newtheorem{theorem}{Theorem}[section]
\newtheorem{algorithm}[theorem]{Theorem/Algorithm}

\newtheorem*{Fconjecture}{$F-$Conjecture}
\newtheorem*{MFconjecture}{$MF-$Conjecture}
\newtheorem{corollary*}{Corollary}
\newtheorem{corollary}[theorem]{Corollary}

\newtheorem{conjecture*}{Conjecture}
\newtheorem*{theorem*}{Theorem}
\newtheorem*{bridgetheorem*}{The Bridge Theorem}
\newtheorem*{raytheorem*}{The Ray Theorem}
\newtheorem*{lemma*}{Lemma}
\newtheorem{lemma}[theorem]{Lemma}

\newtheorem*{claim*}{Claim}
\newtheorem{note}[theorem]{Note}
\newtheorem*{example*}{Example}

\newtheorem*{proof*}{proof}
\newtheorem*{prop*}{Proposition}

\newtheorem{Definition/Proposition}[theorem]{Definition/Proposition}
\newtheorem{definition/lemma}[theorem]{Definition/Lemma}
\newtheorem{definition/theorem}[theorem]{Definition/Theorem}
\newtheorem{example/lemma}[theorem]{Example/Lemma}

\newtheorem*{question*}{Question}

\newtheorem{definition}[theorem]{Definition}

\def\M{{\overline{\mathnormal{M}}}}

\begin{document}

\large
   \begin{center}
       \textbf{Numerical criteria for divisors on $\M_{g}$ to be ample}
      \\[.2em]
      \textbf{}
        \\[1.5em]
       \textbf{Angela Gibney} \\[1.5em]
						        \textbf{} \\[1.5em]
										
   \end{center}

\noindent
\textit{Abstract} The moduli space $\M_{g,n}$ of $n-$pointed stable curves of genus $g$ is stratified by the topological type of 
the curves being parametrized: the closure of the locus of curves with $k$ nodes has codimension $k$.   The 
one dimensional
components of this stratification are smooth rational curves 
(whose numerical equivalence classes are) called  $F-$curves.  These  are believed to determine all ample divisors:
\begin{Fconjecture}A divisor on $\M_{g,n}$ is ample if and only if it positively intersects the $F-$curves.
\end{Fconjecture}
In this paper the $F-$conjecture on $\M_{g,n}$ is reduced to showing that certain divisors in $\M_{0,N}$ for $N \le g+n$ are 
equivalent to the sum of the canonical divisor plus an effective divisor supported on the boundary (cf. Theorem \ref{MFimpliesF}).  
As an application of  the reduction,  numerical criteria are given which if satisfied by a divisor $D$ on $\M_g$, show
that $D$ is ample (cf. Corollaries \ref{$b_0$},\ref{level0criteria}, \ref{$b_1$},  \ref{$b_m$}, and \ref{induct}).  Additionally, an
algorithm is described to check that a given divisor is ample (cf. Theorem/Algorithm \ref{sometimes}).  Using a computer program called  The Nef Wizard, written by Daniel Krashen,  one can use  the criteria and the algorithm to verify the conjecture for low genus.  This
is done on $\M_g$ for $g \le 24$, more than doubling the known cases of the conjecture and showing it is true for the first genus such that $\M_g$ is known to be of general type.

\section{\textbf{Introduction}}
 The moduli
 space $M_{g,n}$ of smooth $n-$pointed curves of genus $g$, and its projective closure, the Deligne-Mumford compactification $\M_{g,n}$ have been
 studied
 in many areas of mathematics.   This is because often properties of families of curves may 
 be translated into facts about the
birational geometry of the moduli space.   For
example,  asking whether almost any curve of genus $g$ occurs as a member of a family given by \emph{free}
parameters -- i.e. parametrized by an open subset of affine space -- is the same as asking whether $\M_{g,0}$ (just written $\M_{g}$) is
unirational.

To learn about the birational geometry of a projective variety like $\M_{g,n}$, it is useful to study its nef and effective
divisors.  A nef divisor $D$ on a projective variety $X$ is a divisor that nonnegatively intersects every effective curve on $X$.
The nef divisors on $X$
 parametrize
 morphisms from $X$ to any projective variety since 
 to every regular map $f: X \longrightarrow Y$ from $X$ to a projective variety $Y$ there corresponds a nef divisor, $f^*A$, where $A$ is
 ample  on $Y$.  The nef and effective divisors of a variety $X$ form cones inside the N\'{e}ron-Severi space of $X$.  Interior to the nef cone is the cone of ample divisors.  By studying these
cones, one can say a lot about the space $X$.  For example, one of the strongest results 
about the birational geometry of $\M_g$ is that for $g \ge 24$, the moduli space is of general type. 
This was proved by Harris, Mumford (and later Eisenbud) who after learning enough about the cone of effective divisors
were able to show that for $g\ge 24$ the canonical divisor of $\M_g$ is interior to it,  and furthermore does not touch the sides.
In particular, in this range $\M_g$ is not unirational, and so the
general curve of genus $g \ge 24$ does not appear as a member of a family of curves parametrized by an open subset
of affine space.

The cone of nef divisors of a projective variety $X$ is always contained inside the effective cone of divisors of $X$. 
For $\M_g$, the nef cone is strictly interior to the effective cone in the sense that they only intersect
at the origin (cf. \cite{thesis}).  As a result of this fact, there is no projective
morphism with connected fiber from $\M_g$ to any lower dimensional variety other than a point.  This is another example given to illustrate 
that the cones of nef and effective divisors 
are extremely important tools for understanding the birational geometry of a projective variety $X$.  Much more information
would be gained if one could further clarify the relationship between the nef and effective cones.  Ideally one would like to 
describe the nef cone explicitly.

One might hope to  specify which
divisors on a projective variety  $X$ are nef by finding a collection of curves $\{C_i\}_{i \in I}$ which 
determine all effective curves -- i.e. span the extremal rays of the Mori Cone of curves.
If such a collection of curves exists,  then one could say that a divisor $D$ on $X$ is nef if and only if it intersects them.
Finding such curves for a given variety $X$ is a very difficult and often 
impossible task.  However for $\M_{g,n}$ there are smooth rational curves called $F-$curves 
that seem to be the right ones to consider.

 In order to describe the $F-$curves, a few facts about the structure of $\M_{g,n}$ will be given. 
  Points in  $\M_{g,n}$ correspond to stable $n-$pointed curves of genus $g$.  A stable curve has at worst nodal singularities.
	 The locus of curves with $k$ nodes has codimension $k$ in $\M_{g,n}$.  Since the dimension of $\M_{g,n}$ is $3g-3+n$, the (closure of the) locus of curves with $3g-4+n$ nodes is $1$-dimensional.  Any curve that is numerically equivalent to a component of this
	  $1$-dimensional locus is called an $F-$curve.    An $F-$divisor is any divisor that nonnegatively intersects all the $F-$curves.  The F-conjecture		 asserts that the F-cone
		  of divisors is the same as the nef cone of divisors of $\M_{g,n}$.

\begin{Fconjecture}A divisor on $\M_{g,n}$ is nef if and only if it nonnegatively intersects a class of curves
called the $F-$curves.
\end{Fconjecture}

In this paper the $F-$conjecture on $\M_{g,n}$ is reduced to showing that certain divisors in $\M_{0,N}$ for $N \le g+n$ are 
equivalent to the sum of the canonical divisor plus an effective divisor supported on the boundary (cf. Theorem \ref{MFimpliesF}).  
As an application of  the reduction,  numerical criteria are given which if satisfied by a divisor $D$ on $\M_g$, show
that $D$ is nef (cf. Corollaries \ref{$b_0$}, \ref{level0criteria}, \ref{$b_1$},  \ref{$b_m$}, and \ref{induct}).  An algorithm is described
for using the reduction to check that a given $F-$divisor is nef (cf. Theorem/Algorithm \ref{sometimes}).  Using a computer program called  The Nef Wizard,  one can show that the criteria and the algorithm completely determine all nef divisors on $\M_g$ for $g \le 24$.  The computer package written by Daniel Krashen can be found at {\textit{http://www.math.yale.edu/users/dkrashen/nefwiz/}}.       

Most of the criteria are phrased so that they can be applied to showing that $F-$divisors on $\M_g$ are nef.  However, since by \cite{GKM}, Theorem $.7$, any $F-$divisor in $\M_{0,g}/S_g$ is the pullback of an $F-$divisor on $\M_g$, they can also be used to prove that $F-$divisors on this space are nef.

It is worth noting that since as it turns out there are a  finite number $F-$curves to begin with, if the $F-$Conjecture is true, then it means that there are finitely many extremal rays of the cone of  curves.  This is  surprising
 since the most general thing one can say about the shape of the cone of curves for an arbitrary  variety 
$X$,  is that the part of the cone corresponding to curves which negatively intersect the canonical divisor is polyhedral; on this part of the cone there are countably many extremal 
 rays and they are spanned by irreducible, rational curves.  The cone of curves for $\M_{g,n}$ is not 
$K-$negative; in fact, since  for $n=0$ and $g \ge 24$, the space $\M_g$ is of general type, 
very much the opposite is true.   It is for this reason that in this work the $F-$conjecture is checked for genus up to $24$;
for higher genera there doesn't seem to be any feature of the spaces which might prevent the conjecture from being true.  Also,
the list of generators of the cone of $F-$ divisors grows extremely fast and so it takes the computer a long time to run through the list of divisors to check
that the criteria are met and the divisors are nef.    

\bigskip

\noindent
{\textit{Previous Results: \ \ }}
Prior to this work, the $F-$conjecture
was  known to be true on $\M_g$ for $g \le 11$ and for 
$g=13$.   
The first cases $g =3$ and $4$ were proved by Carel
Faber for whom the $F-$curves and divisors are 
named.   In \cite{GKM} it is shown that the problem of describing the nef divisors 
on $\M_{g,n}$ can be reduced to solving the $F-$Conjecture on $\M_{0,g+n}$.  
Results of Keel and McKernan \cite{KM} when combined with \cite{GKM} prove
the conjecture for $g \le 11$.  
Farkas and I were able to extend their results to $g = 13$. 

\bigskip

{\textit{Acknowledgements}} I'd like to thank Carel Faber, Bill Fulton, Se\'{a}n Keel, Karen Smith and  Gavril Farkas 
for comments on this work.  I'd like to thank Ravi Vakil for suggesting
using computer to check the criteria.  Daniel Krashen carried out the extremely time consuming and creative effort of 
writing  Nef Wizard.  Our brainstorming about how to design things so that the cases of the $F-$conjecture up through $24$ could be checked in a reasonable amount of time made it happen. 

\section{Definitions and Notation}
Standard definitions are used for cones of  divisors and curves  as well as 
for the basic divisor classes on  $\M_{g,n}$ (cf. eg. \cite{Kollar}, \cite{GKM}, \cite{GF}).  Since numerical details are
 referred to specifically,
 the $F-$curves and divisors
 will now be defined.  Following that, formulae for the  pullback of a divisor
 along certain morphisms will be derived. Since the formulas in Sections $2.2$, $2.3$ and $2.4$ are so  combinatorially involved, one  may 
 wish to skip ahead to Section $3$ and refer back as necessary.

 \subsection{Faber Curves and Divisors}
 
An $F-$curve in $\M_{g,n}$ is any curve that is numerically equivalent to a component of the locus of points in $\M_{g,n}$ having 
$3g-4+n$ nodes.   A subset of the boundary classes $\delta_{i,I}$, taken together with
the tautological classes $\psi_i=-\delta_{i,\emptyset}$ along with the Hodge 
class $\lambda$  
form a basis for the Picard group of $\M_{g,n}$.  By writing
a divisor in terms of these classes and intersecting it with the various $F-$curves, one can see that
if the divisor is an $F-$divisor then its coefficients satisfy certain inequalities. These
 inequalities, which can be taken to define an $F-$divisor, are listed  below.

\begin{definition/theorem}\label{F-inequalities}(cf. \cite{GKM}, Theorems $2.1$ and $2.2$) For $N = \{1 \ldots n\}$
$$D= a \lambda - b_0 \delta_0- \sum_{
\stackrel{0 \le i \le \lfloor \frac{g}{2} \rfloor}{\stackrel{ I \subseteq N}{\text{if } i=0, |I| \ge 1}}
}b_{i,I} \delta_{i,I}$$
is an $F-$divisor on $\M_{g,n}$ if and only if:
\begin{enumerate}
\item $a-12b_0+b_{1,\emptyset}$,
\item $b_{i,I} \ge 0$,
\item $2b_0 - b_{i,I} \ge 0$,
\item $b_{i,I} + b_{j,J} \ge b_{i+j, I \cup J},$ for all $i$ and $j$ such that $i+j \le n-1$, and such that $I \cap J = \emptyset$,
\item $b_{i,I} + b_{j,J} + b_{k,K} + b_{l,L} - ( b_{i+j, I \cup J} +  b_{i+k, I \cup L}+
b_{i+l, I \cup L}) \ge 0$, for all $i$, $j$, $k$ and $l$ such that $i+j+k+l=g$, and $I \cup J \cup K \cup L = \{1\ldots n\}$.
\end{enumerate}
\end{definition/theorem} 

\subsection{Boundary Restrictions}

\noindent
Let $f:\M_{0,g+n}\longrightarrow \M_{g,n}$ be the morphism obtained by attaching a pointed curve of genus $1$ to
each of the first $g$ marked points.
The pullback  $f^*D$ will often be referred to as the restriction of a divisor $D$ to the flag locus.

Divisors in $\M_{0,g}$ pulled back along certain
so-called boundary restriction morphisms will also be considered.  

\begin{definition}\label{bdryrest}For $z \ge 2$, $i \ge 1$ and disjoint subsets $N_j \subset N = \{1 \ldots n \}$,
of order $n_j \ge 2$, let $[N_1:N_2:\ldots :N_a]$ be the boundary restriction morphism which we denote by
$v_{a,z}:\M_{0,a+z} \longrightarrow \M_{0,n}$, where $n=\sum_{j=1}^{a}n_j + z$ given by attaching an $(n_j +1)-$pointed
genus $0$ 
curve (whose marked points consist of an attaching point and the  $N_j$) to each of the first {\textbf{a}} marked points and
doing nothing to the last {\textbf{z}} marked points.  We say that \textbf{$a$} is the order of the boundary restriction morphism.  
\end{definition}

Note  that if $D$ is an $F-$divisor in $\M_g$ then $f^*D$ is an F divisor in $\M_{0,g}$.  
Likewise, if $D$ is an F divisor in $\M_{0,n}$ and 
$v: \M_{0,a+z} \longrightarrow \M_{0,n}$ is a boundary restriction, then $v^*D$ is an $F-$divisor in $\M_{0,a+z}$. 
As shown in Lemma $1$, for any boundary restriction morphism $v=[N_1 \ldots N_a]$, 
the pullback $v^*f^*D$ is determined by the orders of the sets $N_j$.  Hence one may denote the boundary restriction
morphism $[N_1 \ldots N_a]$ by the $a-$tuple $[n_1 \ldots n_a]$.

\begin{lemma}\label{pullback}
\begin{enumerate}
\item Let $D = a \lambda -\sum_{i=0}^{\lfloor \frac{g}{2} \rfloor}b_i \delta_i$ in $\M_g$ be a divisor.  Then $$f^*D = b_1 \sum_{i=1}^{g}\psi_i - \sum_{i=2}^{\lfloor \frac{g}{2} \rfloor}b_i B_i,  \mbox{\ \ \ \ where  \ \  } B_i = \sum_{\stackrel{I \subset \{1\ldots g\}}{ |I|=i}}\delta_{I}.$$
\item Let
 $v=v_{a,z}=[n_1 \ldots n_a]$ be a boundary restriction of $\M_{0,g}$.  
 Then $$v^*f^*D=\sum_{i\in A=\{i|n_i \ge 2\}}b_{n_i}\psi_{i}+b_1\sum_{i\in Z=A^c}\psi_{i}-\sum_{\stackrel{B \subset A, 0 \le y \le |Z|}{2\le y+|B|\le \lfloor \frac{a+z}{2}\rfloor}} b_{y+\sum_{i\in B}n_i} \Delta^{Z,y}_{A,B},$$ where $\Delta^{Z,y}_{A,B}=\sum_{Y \subset  Z, |Y|=y}^{}\delta_{Y \cup B}$.
 \end{enumerate}
\end{lemma}

\begin{proof}
For both formulas apply Lemma $1.4$ in \cite{AC}, page $5$. 
\end{proof}

\subsection{Description of the $A-$averages}

The equivalence classes of boundary divisors span  $\mbox{Pic}(\M_{0,n})$  but are not independent.  
Consequently, any divisor class in $\M_{0,n}$ such as the  $\psi_i$   
can be expressed in terms of the boundary classes and moreover, there are different ways of doing so.
Given $i$, $j$ and $k \in \{1\ldots n\}$,  one has that $$\psi_i = \sum_{\stackrel{I \subseteq \{1 \ldots n\}}{i \in I; j,k \notin I}}\delta_I.$$
\noindent
In particular, there are $\binom{n-1}{2}$ ways of expressing a divisor class $\psi_i$ as a sum of boundary divisors
in this manner.  By combining these
in various ways one can produce different manifestations of the $\psi_i$ as sums of boundary classes.  Suppose $A \subseteq \{1 \ldots n \}$, and $i \in A$.  In this section, four ways to write  $\psi_i$
 as a sum of boundary divisors with respect to $A$ will be given.  
These are used to express a general divisor $D$ on $\M_{0,n}$ in terms of boundary classes and 
enable one to locate where the divisor sits in the N\'{e}ron-Severi space of $\M_{0,n}$ with respect
to its effective cone of divisors.

\begin{note}For $B \subset A$, we will often use the notation $$\Delta^{Z,y}_{A,B}=\sum_{Y \subset Z, |Y|=y}\delta_{Y \cup B}.$$
\end{note}

\noindent
The first way to write  $\psi_i$, for $i \in A$ as a sum of boundary classes comes from combining all the expressions for $\psi_i$ given
above such that $j$,$k \in A \setminus \{i\}$.

\begin{definition/lemma} \label{A1} Let $A \subseteq \{1 \ldots n\}$ with $a=|A| \ge 3$ and $Z = A^c$ with $z=|Z|$.  The 
{\textbf{1-st A-average of $\psi_i$ with $i \in A$}} is 
$$\psi_i = \sum_{\stackrel{i \in B \subset A}{|B|=b}}\frac{(a-b)(a-b-1)}{(a-1)(a-2)}\Delta^{Z,y}_{A,B}.$$
\end{definition/lemma}

\noindent
The second $A-$average is derived by writing down all such expressions for $\psi_i$, by taking pairs $j \in A \setminus \{i\}$ and
$k \in A^c$.

\begin{definition/lemma} \label{A2} Let $A \subset \{1 \ldots n\}$ with $a=|A| \ge 2$ and $Z = A^c$ with $z=|Z|\ge 1$.  The 
{\textbf{2nd A-average of $\psi_i$ with $i \in A$}} is 
$$\psi_i = \sum_{\stackrel{i \in B \subset A}{|B|=b}}\frac{(a-b)(z-y)}{(a-1)z}\Delta^{Z,y}_{A,B}.$$
\end{definition/lemma}

\noindent
The third $A-$average of $\psi_i$ is generated by taking all expressions such that $j$,$k \in A^c$.

\begin{definition/lemma} \label{A3} Let $A \subset \{1 \ldots n\}$ with $a=|A| \ge 1$ and $Z = A^c$ with $z=|Z|\ge 2$.  The 
{\textbf{3rd A-average of $\psi_i$ with $i \in A$}} is 
$$\psi_i = \sum_{\stackrel{i \in B \subset A}{|B|=b}}\frac{(z-y)(z-y-1)}{z(z-1)}\Delta^{Z,y}_{A,B}.$$
\end{definition/lemma}

\noindent
Finally, by taking all possible pairs $j$ and $k \in \{1 \ldots n \}\setminus \{i\}$, one obtains the fourth
expression for $\psi_i$ in terms of the boundary classes.  This  $A-$average comes from taking the largest number of ways of expressing the $\psi_i$ as a sum of boundary divisors in this
way and is referred to as the {\textbf{big average of $\psi_i$}}.  It can be found in \cite{GF}, Lemma $1$.

\begin{definition/lemma} \label{A4} \label{bigA} For $i \subset \{1 \ldots n\}$, the 
{\textbf{4th or big A-average of $\psi_i$}} is 
$$\psi_i = \sum_{Y \subset \{1\ldots n\}\setminus \{i\}}\frac{(n-1-y)(n-2-y)}{(n-1)(n-2)}\delta_{Y \cup \{i\}},$$
where $y$ is the number of elements in the set $Y$.
\end{definition/lemma}

\subsection{The $c-$averages of a divisor on $\M_{0,n}$}
The main technique in this work is to use different ways to write certain divisors on $\M_{0,n}$ as $$cK_{\M_{0,n}}+E,$$
where $E$ is an effective sum of boundary classes.  These expressions are called $c-$averages of a divisor.

For
$D=a\lambda-\sum_{2 \le i \le \lfloor \frac{g}{2} \rfloor }b_i \delta_i,$ on $\M_g$, the pull-back $v^*f^*D$ on $\M_{0,n}$
of $D$ along the so called boundary restriction morphisms
can be expressed as 

 $$v^*f^*D=b_1\sum_{i\in Z=\{i|n_i=1\}}\psi_{i}+\sum_{i\in A=\{i|n_i \ge 2\}}b_{n_i}\psi_{i}-\sum_{\stackrel{B \subset A, 0 \le y \le |Z|}{2\le y+|B|\le \lfloor \frac{a+z}{2}\rfloor}} b_{y+\sum_{i\in B}n_i} \Delta^{Z,y}_{A,B}.$$ 

 By replacing the $\psi_i$ in the expression above for $i \in A$ (respectively for $i \in Z$) with combinations
 of the various $A-$
and  $Z-$averages one obtains up to $12$ different $c-$averages of the divisor $v^*f^*D$.  When $a=0$, there is just the big average.   To give a flavor for what the expressions
look like, three examples are given below.





\begin{definition/lemma} Suppose $D$ is a divisor in $\M_g$ and $c \ge 0$. Let  $f:\M_{0,g} \longrightarrow \M_g$ be the morphism given
by attaching elliptic tails and $v:\M_{0,a+z} \longrightarrow \M_{0,g}$ be a boundary restriction morphism.
 The  {\bf{big c-average}} of $v^*f^*D$ is:
 \begin{multline*}
 v^*f^*D=c K_{\M_{0,a+z}}
 + \sum_{\stackrel{B \subseteq A, |B|=b}{ 0 \le y \le z, 2 \le y+b \le \lfloor \frac{a+z}{2} \rfloor}} \big(f_{y,b} (b_1-c)+ 
 g_{y,b}  \sum_{i \in B} (b_{n_i}-c) 
 \\
 + h_{y,b} \sum_{i \in B^c} (b_{n_i} -c) + 2c - b_{y+\sum_{i \in B}b_{n_i}}\big) \Delta^{Z,y}_{A,B},
 \end{multline*}
 where $g_{y,b}=\frac{(a+z-y-b)(a+z-y-b-1)}{(a+z-1)(a+z-2)}$, $h_{y,b}=\frac{(y+b)(y+b-1)}{(a+z-1)(a+z-2)}$,  $f_{y,b} = y  g_{y,b} + (z-y) h_{y,b}$ and $\Delta^{Z,y}_{A,B}=\sum_{Y \subset Z, |Y|=y}\delta_{Y \cup B}$.
 \end{definition/lemma}

 \begin{proof}

Recall that from Lemma $1$, if
 $v=v_{a,z}=[n_1 \ldots n_a]$ is a boundary restriction of $\M_{0,g}$, then $$v^*f^*D=b_1\sum_{i\in Z=\{i|n_i=1\}}\psi_{i}+\sum_{i\in A=\{i|n_i \ge 2\}}b_{n_i}\psi_{i}-\sum_{\stackrel{B \subset A, 0 \le y \le |Z|}{2\le y+|B|\le \lfloor \frac{a+z}{2}\rfloor}} b_{y+\sum_{i\in B}n_i} \Delta^{Z,y}_{A,B}.$$ 
Using the relation  $K_{\M_{0,g}}=\sum_{1 \le i \le g}^{}\psi_i-2 \Delta$,  rewrite the expression as:
\begin{multline*}
v^*f^*D=\sum_{j \in A}^{}(b_{n_j}-c)\psi_j+(b_1-c)\sum_{j \in Z}^{}\psi_j
+ c \sum_{j \in A \cup Z}^{}\psi_j-\sum_{\stackrel{0\le y \le z, B \subseteq A}{2 \le y+|B| \le \lfloor \frac{n}{2} \rfloor}}^{}b_{y+\sum_{k\in B}n_k}  \Delta^{Z,y}_{A,B}\\
=\sum_{j \in A}^{}(b_{n_j}-c)\psi_j+(b_1-c)\sum_{j \in Z}^{}\psi_j
+ c K_{\M_{0,g}}+  \sum_{\stackrel{0\le y \le z, B \subseteq A}{2 \le y+|B| \le \lfloor \frac{n}{2} \rfloor}}^{} 
\big( 2c-b_{y+\sum_{k\in B}n_k} \big) \Delta^{Z,y}_{A,B}.
\end{multline*}

By big averaging  the $\psi_i$, and distributing the coefficients through the sum, one obtains the expression given in the
Theorem.
\end{proof}

Note that if $a=0$, then $z=g$, and  $v^*f^*D=f^*D$.  Therefore the  {\bf{big c-average of $f^*D$}}  is just:
 $$f^*D  =c K_{\M_{0,n}}+ \sum_{2 \le i \le \lfloor \frac{g}{2}\rfloor} \big( (b_1-c)\frac{i(g-i)}{(g-1)} + 2c - b_i  \big) B_i.$$  Moreover, as long as $c >0$ one can write this as:
 $$f^*D = c\Big( K_{\M_{0,n}}+ \sum_{2 \le i \le \lfloor \frac{g}{2}\rfloor} \big( \frac{2g-2-i(g-i)}{g-1}+\frac{b_1i(g-i)-b_i(g-1)}{(g-1)c} \big) B_i \Big).$$

\begin{definition/lemma} Suppose $D$ is a divisor in $\M_g$ and $c \ge 0$. Let  $f:\M_{0,g} \longrightarrow \M_g$ be the morphism given
by attaching elliptic tails and $v=[n_1 \ldots n_a]:\M_{0,a+z} \longrightarrow \M_{0,g}$ be a boundary restriction morphism such that $a \ge 2$
and $z \ge 2$.
 The  {\bf{second  c-average}} of $v^*f^*D$ is:
  
$$v^*f^*D=c K_{\M_{0,a+z}}
	  + \sum_{\stackrel{B \subseteq A, |B|=b}{ 0 \le y \le z, 2 \le y+b \le \lfloor \frac{a+z}{2} \rfloor}} C_{y,B} \ \ \ \Delta^{Z,y}_{A,B},$$
where 
\begin{multline*}
C_{y,B} = \frac{(a-b)(z-y)\sum_{i \in B}(b_{n_i}-c)+ b y \sum_{i \in B^c}(b_{n_i}-c)}{(a-1)z}\\
+ \frac{(b_1-c)\big((a-b)(z-y)+ b y \big)}{a(z-1)} +2c-b_{y+ \sum_{i \in B}n_i}.
\end{multline*}
\end{definition/lemma}

					 \begin{proof}

 Recall that from Lemma $1$, if
					  $v=v_{a,z}=[n_1 \ldots n_a]$ is a boundary restriction of $\M_{0,g}$, then $$v^*f^*D=b_1\sum_{i\in Z=\{i|n_i=1\}}\psi_{i}+\sum_{i\in A=\{i|n_i \ge 2\}}b_{n_i}\psi_{i}-\sum_{\stackrel{B \subset A, 0 \le y \le |Z|}{2\le y+|B|\le \lfloor \frac{a+z}{2}\rfloor}} b_{y+\sum_{i\in B}n_i} \Delta^{Z,y}_{A,B},$$ where $\Delta^{Z,y}_{A,B}=\sum_{Y \subset  Z, |Y|=y}^{} \delta_{Y \cup B}$.
						 
	Using the relation  $K_{\M_{0,g}}=\sum_{1 \le i \le g}^{}\psi_i-2 \Delta$,  rewrite the expression as:
						 \begin{multline*}
						 v^*f^*D=\sum_{j \in A}^{}(b_{n_j}-c)\psi_j+(b_1-c)\sum_{j \in Z}^{}\psi_j
						 + c \sum_{j \in A \cup Z}^{}\psi_j-\sum_{\stackrel{0\le y \le z, B \subseteq A}{2 \le y+|B| \le \lfloor \frac{n}{2} \rfloor}}^{}b_{y+\sum_{k\in B}n_k}  \Delta^{Z,y}_{A,B}\\
						 =\sum_{j \in A}^{}(b_{n_j}-c)\psi_j+(b_1-c)\sum_{j \in Z}^{}\psi_j
						 + c K_{\M_{0,g}}+  \sum_{\stackrel{0\le y \le z, B \subseteq A}{2 \le y+|B| \le \lfloor \frac{n}{2} \rfloor}}^{} 
						 \big( 2c-b_{y+\sum_{k\in B}n_k} \big) \Delta^{Z,y}_{A,B}.
						 \end{multline*}

						 By substituting the second $A-$average of the $\psi_i$ for $i \in A$, with $Z = A^c$ and distributing the coefficients through the sum and by substituting the second $Z-$average of the $\psi_i$ for $i \in Z$, with $A = Z^c$ and distributing the coefficients through the sum, one obtains the expression given in the
						 Theorem.
						 \end{proof}

One could also combine different averages of the $\psi_i$.  For example by substituting the second $A-$average of the $\psi_i$ for $i \in A$, and by substituting the big $Z-$average of the $\psi_i$ for $i \in Z$, with $A = Z^c$ and distributing the coefficients through the sum, one obtains the expression:

$$v^*f^*D=c\big( K_{\M_{0,a+z}}
	  + \sum_{\stackrel{B \subseteq A, |B|=b}{ 0 \le y \le z, 2 \le y+b \le \lfloor \frac{a+z}{2} \rfloor}} 
		\frac{n(\alpha)}{d(\alpha)} + \frac{n(\beta)}{d(\beta)}/c \ \ \ \Delta^{Z,y}_{A,B}\big),$$
where 
\begin{multline*}
n(\alpha) = (a+z-1)(a+z-2)(2(a-1)-b(a-b))\\
-(a-1)((a+z-y-b)(a+z-y-b-1)y+(y+b)(y+b-1)(z-y)),
\end{multline*}
$$d(\alpha)=(a-1)(a+z-1)(a+z-2),$$
\begin{multline*}
n(\beta) = (a+z-1)(a+z-2)((z-y)(a-b)\sum_{i \in B}b_{n_i}+y \ b\sum_{i \in B^c}b_{n_i})\\
+b_1(a-1)z((a+z-y-b)(a+z-y-b-1)y+(y+b)(y+b-1)(z-y)).
\end{multline*}
and $$d(\beta) = (a+z-1)(a+z-2)(a-1)z.$$

\section{Reduction of the $F-$Conjecture}
\noindent
In this section the $F-$Conjecture is reduced to what will be termed the $MF-$Conjecture  which asserts that $F-$divisors on $\M_{0,N}$
are the sum of the canonical divisor and an effective divisor.
\begin{MFconjecture}
Every $F-$divisor on $\M_{0,N}$ is of the form $cK_{\M_{0,N}}+E$ where $c \ge 0$ and $E$ is an effective sum of boundary classes.
\end{MFconjecture}

\noindent
The problem of whether every $F-$divisor in $\M_{0,N}$ is numerically equivalent to an effective sum of boundary classes
has been called Fulton's Conjecture (cf. \cite{GF}).  Here the $MF$ stands for "Modified Fulton's" named so because 
Fulton's Conjecture is  the  $MF-$Conjecture  with $c=0$.  

\medskip
\noindent
The numerical criteria and 
 algorithm for showing a divisor is nef given in the next section rest on this reduction of the $F-$conjecture to the $MF-$conjecture. 
 \begin{theorem}\label{MFimpliesF}
If the $MF-$Conjecture is true on $\M_{0,N}$ for $N \le g+n$, then the $F-$conjecture is true on $\M_{g,n}$.  In particular,
if the $MF-$conjecture is true then the $F-$conjecture is true.
\end{theorem}

\noindent
Two facts are needed to explain how the $MF-$Conjecture implies the $F-$conjecture.  The first is that if the $F-$conjecture is true in $\M_{0,g+n}$ then it is true in $\M_{g,n}$.
More precisely, let $f: \M_{0,g+n} \longrightarrow \M_{g,n}$ be the morphism associated to the map given by attaching pointed elliptic tails at each of the first $g$ marked points. 

\begin{bridgetheorem*}(\cite{GKM}, Thm $.03$) 
A divisor $D$ on $\M_{g,n}$ is nef if and only if
\begin{enumerate}
\item $D$ is an $F-$divisor, and 
\item $f^*D$ is a nef divisor on $\M_{0,g+n}$.
\end{enumerate}
\end{bridgetheorem*}

\noindent
The following result is the second important fact needed to prove Theorem \ref{MFimpliesF}.
\begin{raytheorem*} (\cite{GF}, Thm $4$ and \cite{KM}, Thm $1.2$) If $R$ is an extremal ray of the cone of curves of $\M_{0,N}$ and if 
$R \cdot  (K_{\M_{0,N}}+G) < 0$ where $G$ is any effective sum of boundary
for which $\Delta \setminus G$ is nonnegative, then $R$ is spanned by an F curve.
\end{raytheorem*}

\noindent
The symbol $\Delta$ denotes the sum of boundary classes.  So the condition in the Ray Theorem  is that $G = \sum_{S}a_S \delta_S$ such that  $0 \le a_S \le 1$ for all $S$.  The Ray Theorem  is an 
extension of work by  Keel and McKernan 
which states that if $R$ is an extremal ray of $\overline{NE}(\M_{0,N})$
and if  $R \cdot  (K_{\M_{0,N}}+G) \le 0$ for 
$G = \sum_{S}a_S \delta_S$ such that  $0 \le a_S < 1$, then $R$ is spanned by an F curve.  

\bigskip
\begin{proof}(of Theorem \ref{MFimpliesF}) \ Suppose that whenever one has an $F-$divisor $D$ on $\M_{0,N}$,   there exists a constant $c \ge 0$ for which 
$$D = cK_{\M_{0,N}}+E,$$
where $E$ is an effective sum of boundary classes.  We will show that this assumption implies that the $F-$conjecture is true on
$\M_{g,n}$.  By the Bridge Theorem, in order to prove the $F-$conjecture on $\M_{g,n}$, it is enough to show that any $F-$divisor
on $\M_{0,g+n}$ is nef. Hence if we  show that our assumption implies that $D$ is nef, then the theorem is proved.

By definition, if $D$ nonnegatively intersects all the extremal rays of the cone of curves, then $D$ is nef.  Suppose $R$ is an extremal ray of the cone of curves.  The first thing to note is that since $D$ is an $F-$divisor,
and if $R$ is spanned by an $F-$curve, then $D$ nonnegatively intersects $R$.  We will prove that there are no other kinds of extremal rays.  We do this by induction on the number of marked points. 
As base case we take $N=7$ since the $F-$conjecture is true for $N \le 7$ (cf. \cite{KM}).

The cone of curves is the closure of NE$(\M_{0,N})$ in
the real vector space $N_1(\M_{0,N})$.  So every extremal ray $R$ is either spanned by an irreducible curve or is the limit of rays spanned by irreducible curves.

Suppose that $R$ is a $D-$negative extremal ray of the cone of curves of $\M_{0,N}$ for $N > 7$ that isn't spanned by an $F-$curve.
In other words, suppose that  $$R \cdot D = R \cdot ( c K_{\M_{0,N}} + E) < 0.$$
In particular, by The Ray Theorem, $ R \cdot  E < 0$.

If $R$ is spanned by a curve, then since $E$ is an effective sum of boundary classes,
to get a contradiction, 
it is enough to show that $D$ is nef when restricted to the components in the support of $E$.  This results in pulling $D$ back
to a space $\M_{0,n}$, for $n < N$ along a boundary restriction morphism (defined in Section $2$).  Since the pullback of
an $F-$divisor along a boundary restriction morphism is an $F-$divisor, one can repeat this argument until ending up in $\M_{0,n}$ for $n \le 7$.

If the extremal ray $R$ is a limit of curves, then one can find a ray $R'$ spanned by a curve that is close enough so that $R'$ intersects $D$ negatively.  In this case one reaches a contradiction as above.
\end{proof}

As is shown in Theorem \ref{Fulton} below, the $MF-$Conjecture  is true on $\M_{0,N}$ for $N \le 6$.   It was already known
to be true with $c=0$ ( cf. \cite{GF}).  However,
 the proof with $c=0$ is much harder since showing that a divisor class is in the convex hull of boundary classes is more difficult than showing it is in the convex hull of boundary classes and the canonical divisor.  It seems unlikely that $MF$ is true with $c=0$, even
 when $N=7$.

\begin{theorem}\label{Fulton}
If $D$ is any divisor on $\M_{0,n}$ for $n=5$ or $6$, then there exists a constant $k >0$ such that $D = cK_{\M_{0,n}}+E$,
for all $c \ge k$ and where $E$ is an effective sum of boundary divisors.  In particular, the $MF-$conjecture is true on 
$\M_{0,n}$ for $n \le 6$.
\end{theorem}

			\begin{proof}(of Theorem \ref{Fulton}) \   First suppose that $n=5$.  By substituting the big averages (see Section $2$ for definitions)
			of the divisors $\psi_i$ one can express the divisor $D$ as follows:
			\begin{multline*}
			D=\sum_{1 \le i \le 5}c_i \psi_i = c(\sum_{1 \le i \le 5}\psi_i-2\sum_{ij \in \{1 \ldots 5 \}}\delta_{ij})
			+\sum_{1 \le i \le 5}(c_i-c) \psi_i + 2c \sum_{ij \in \{1 \ldots 5 \}}\delta_{ij}\\
			=cK_{\M_{0,5}}+\sum_{ij \in \{1 \ldots 5\}}\big(\frac{1}{2}\sum_{k \in \{i,j\}}(c_k-c)
			+\frac{1}{6}\sum_{k \in \{i,j\}^c}(c_k-c)+2c  \big)\delta_{ij}\\
			=cK_{\M_{0,5}}+\sum_{ij \in \{1 \ldots 5\}}\big(\frac{1}{2}\sum_{k \in \{i,j\}}c_k
			+\frac{1}{6}\sum_{k \in \{i,j\}^c}c_k+\frac{1}{2}c  \big)\delta_{ij}.
			\end{multline*}
			Similarly, when $n=6$ one can write $D$ as follows:
			\begin{multline*}D=\sum_{1 \le i \le 6}c_i \psi_i-\sum_{ij \in \{2 \ldots 6\}}b_{1ij}\delta_{1ij} \\
			=cK_{\M_{0,6}}+\sum_{ij \in \{1 \ldots 6\}}\big(\frac{1}{2}\sum_{k \in \{i,j\}}(c_k-c)
			+\frac{1}{10}\sum_{1 \le k \le 6}(c_k-c)+2c  \big)\delta_{ij}\\
			+\sum_{ij \in \{2 \ldots 6\}}\big(
			+\frac{3}{10}\sum_{1 \le k \le 6}(c_k-c)+2c-b_{1ij}  \big)\delta_{1ij}\\
			=cK_{\M_{0,6}}+\sum_{ij \in \{1 \ldots 6\}}\big(\frac{1}{2}(c_i+c_j)
			+\frac{1}{10}\sum_{1 \le k \le 6}c_k+\frac{2c}{5} \big)\delta_{ij} \\
			+\sum_{ij \in \{2 \ldots 6\}}\big(
			\frac{3}{10}\sum_{1 \le k \le 6}c_k+\frac{c}{5}-b_{1ij}  \big)\delta_{ij}.\\
			\end{multline*}
			In either case, if $c$ is taken to be big enough, then the assertion is true.  
			\end{proof}
			\bigskip

\section{ Iterative Procedures to Show a Divisor in $\M_g$ is Nef}
By proving particular cases of the $MF-$conjecture, one can use Theorem \ref{MFimpliesF} to define an algorithm
for proving that a divisor $D$ in $\M_{g}$ is nef (cf. Theorem  \ref{sometimesrestrict}).  The first step 
is the following result.

\begin{theorem}\label{step0}
  If $D = b_1\sum_{1 \le i \le g}\psi_i-\sum_{2 \le i \le \lfloor \frac{g}{2}\rfloor}b_i B_i$ is any F divisor in $\M_{0,g}/S_g$, 
	 then there exists a constant $c > 0$ for which $D = cK_{\M_{0,n}}+E$ such that $E$ is an effective sum of boundary classes.
\end{theorem}

\begin{proof}Assume that  $D = b_1\sum_{1 \le i \le g}\psi_i-\sum_{2 \le i \le \lfloor \frac{g}{2}\rfloor}b_i B_i$ is any F divisor in $\M_{0,g}/S_g$ and consider the $c-$average of $D$:
	\begin{multline*}
	D =c\Big(K_{\M_{0,n}} 
	+\sum_{2 \le i \le \lfloor \frac{g}{2}\rfloor}\big( \frac{2g-2-i(g-i)}{g-1}
	  +\frac{b_1i(g-i)-b_i(g-1)}{(g-1)c} \big) B_i \Big).
	\end{multline*}
First assume that $g \ge 8$ and let $c_i = \alpha_i +\frac{\beta_i}{c}$,
where $\alpha_i=\frac{2g-2-i(g-i)}{g-1}$ and $\beta_i=\frac{i(g-i)b_1-(g-1)b_i}{(g-1)}$.
Note that $i(g-i)b_1 > (g-1)b_i$ since $D$ is an $F-$divisor (\cite{KM}, Lemma ?).  In particular, 
$\beta_i > 0$ for all $i$.  One has that  $\alpha_2 = \frac{2}{g-1}$ is positive,
and as shall be shown, $\alpha_i < 0$ for all $i \ge 3$.  Indeed,
the numerator $n(i)=2g-2-i(g-i)$ is negative:
$n(i)$ is decreasing with respect to $i$ since $\frac{\partial n}{\partial i}=i-g$,
and $n(3)=7-g < 0$ for $g \ge 8$.  
																																																																																															 
																																																																		 It shall be argued that there is a positive $c$ for which all the coeficients $c_i$ are positive.  For all $i$ the functions $\alpha_i + \beta_i/c$ have vertical asymptotes at $c=0$. 
																																																																		The function $\alpha_2 + \beta_2/c$ tends to $2/(g-1)$ as $c$ tends to infinity.  For $i > 2$,
 the function $\alpha_i + \beta_i/c$ is concave up and decreasing, crossing the $c$ axis when $c = \beta_i / \alpha_i >0$.
 Hence we can take $$c = \mbox{min}\{\frac{\beta_i}{\alpha_i} | \ \  3 \le i \le \lfloor \frac{g}{2} \rfloor \}.$$  

 In the case $g=7$, one has the relation that $2b_1 \ge b_3$ from intersecting $D$ with the $F-$curve given by the $4-$tuple
$[1:1:2:3]$.  In this case $$D= c\big(K_{\M_{0,7}}+(\frac{c+5b_1-2b_2}{3c})B_2+(\frac{2b_1-b_3}{c})B_3\big).$$
																																																																		In particular, one  must take $c$ so that $c+5b_1-2b_2 \ge 0$.

 Since by Theorem \ref{Fulton}, the result is true more generally for $g \le 6$, the theorem is proved.
																																																																		 \end{proof}

																																																																		 This result was known to be true for $c=0$ (cf. \cite{GF}).  As was pointed out in \cite{GF},
																																																																		 the problem of showing that a particular $F-$divisor $D$ on $\M_g$  is nef can therefore be reduced
	to showing $f^*D=E$ is 	 nef  when restricted to all of the boundary divisors in the support of $E$. 
																																																																			However, as is shown in Theorem  \ref{sometimesrestrict}, that it works for $c > 0$ is	a drastic improvement since one can immediately reduce the problem of showing a particular $F-$divisor is nef to showing it is nef  when restricted to the boundary divisors in the support of $E$ having coefficient larger than $c$.

\begin{theorem}\label{alwaysrestrict} Consider an F-divisor of the form $D=b_1 \sum_{i=1}^{g}\psi_i-\sum_{i=2}^{\lfloor \frac{g}{2} \rfloor}b_iB_i$ in $\M_{0,g}/S_g$.  If  for each boundary restriction $v_{a,z}=[n_1 \dots n_a]$,  there exists a constant $c_v \ge 0$ such that
$$v^*D=c_vK_{\M_{0,a+z}}+E,$$
where $E$ is an effective sum of boundary classes, then $D$ is nef.
\end{theorem}

To prove the next theorem it will be necessary to refer to boundary restriction morphisms
and $c-$averages  which are defined in Section $2.4$.  
In particular, so-called necessary boundary restriction morphisms will be considered.

\begin{definition} Let $D$ be a divisor on $\M_{0,g}$ and suppose that
the $c-$average of $v^*D$
is  of the form 
$cK_{\M_{0,a+z}}+E$, where $c \ge 0$, $E$ is an effective sum of distinct boundary classes, and
  $v_{a,z}=[n_1 \dots n_{a-1}]:\M_{0,a+z} \longrightarrow \M_{0,g}$ is any boundary restriction morphism.  Then
define \textbf{necessary boundary restrictions} to be the  boundary restrictions $v_{S}$ and $v_{S^c}$ such that the coefficient
 of $\delta_S$ in this expression is greater than $c$.  
 Here, for $S \subset \{p_1 \ldots p_{a+z}\}$, one defines $v_S$ to be the boundary restriction
  morphism
	 $$v_S=[\sum_{i \in S \cap A}n_i+|S \cap Z|,\{n_i\}_{i\in S^c \cap A}].$$  Recall that
	 $A=\{p_i \in \{p_1 \ldots p_{a+z}\} | n_i \ge 2 \}$ is the set of attaching points
	 of the boundary restriction morphism and $Z=  \{p_i \in \{p_1 \ldots p_{a+z}\} | n_i =1 \}$ is the set of points to which nothing
	 is attached.  
\end{definition}

\begin{theorem} \label{sometimesrestrict} 
Consider an F-divisor of the form $D=b_1 \sum_{i=1}^{g}\psi_i-\sum_{i=2}^{\lfloor \frac{g}{2} \rfloor}b_iB_i$
in $\M_{0,g}/S_g$.  If for each composition of necessary boundary restrictions $v$,  
there exists a constant $c_v \ge 0$ such that
$$v^*D=c_vK+E,$$
where $E$ is an effective sum of boundary classes, then $D$ is nef.
\end{theorem}

\begin{proof}
By Theorem \ref{step0}, the divisor $D$ is of the form $cK_{\M_{0,g}}+E$ where $c \ge 0$ and $E$ is an effective sum of 
distinct boundary classes.  Exactly as in the proof of Theorem \ref{MFimpliesF}, using The Ray Theorem, 
$D$ is nef as long as it nonnegatively
intersects all curves in the support of any component of $E$ with coefficient larger than $c$.  In other words, suppose the
coefficient of $\delta_S$ is larger $c$, then it is enough 
to show that $D$ is nef when restricted
to both $\Delta_S$.  That is, it is enough to show that $v^*D$ is nef for $v_S=[S]$ and $v_{S^c}=[S^c]$.
By hypothesis, $$v_S^*D=c_vK_{\M_{0,1+g-|S|}}+E,$$
where $E$ is an effective sum of distinct boundary classes. Repeating this argument, it is enough to show that
for each composition of necessary boundary restrictions $v$,  
there exists a constant $c_v \ge 0$ such that
$$v^*D=c_vK+E,$$
where $E$ is an effective sum of boundary classes.  
Eventually the process will stop since the $F-$conjecture is known to be true on $\M_{0,N}$, for $N \le 7$ (cf. \cite{KM}).
\end{proof}

To have a computer check that any composition of necessary boundary restrictions of an $F-$divisor on $\M_{0,g}/S_g$
always restricts to a divisor on $\M_{0,a+z}$ of the form $cK_{\M_{0,a+z}}+E$,
one can use any of the $c-$averages defined in Section $2.5$.  

\begin{algorithm}\label{sometimes} Let $D$ be an F-divisor of the form $a\lambda- \sum_{i=0}^{\lfloor \frac{g}{2} \rfloor}b_i\delta_i$
on $\M_{g}$.  If the $c$ average $v^*f^*D=cK+E$ of any necessary boundary restriction $v$ of $f^*D$ 
is effective, then $D$ is nef.
\end{algorithm}

\begin{proof}
This follows from Theorem \ref{sometimesrestrict}
\end{proof}

\section{Numerical Criteria}
\noindent
In this section,  as an application of Theorem \ref{MFimpliesF} and the iterative procedures given in Section $4$, numerical 
criteria are given which guarantee that divisors on $\M_{g}$ are nef.  These criteria can be viewed as a way of carving the cone of $F-$divisors on $\M_g$ into nef subcones.  As is explained in the following section, these subcones cover the entire $F-$cone
for $g \le 24$.

\begin{corollary}\label{$b_0$}Let $D= a\lambda-\sum_{0\le i \le \lfloor \frac{g}{2} \rfloor}b_i \delta_i$ be
an $F-$divisor on $\M_{g}$.  If  for  $i \in \{ 2, \ldots , \lfloor \frac{g}{2} \rfloor \}$,
      $$-b_0(g-1) \le i(g-i)(b_1-b_0)+(g-1)(b_0-b_i) \le 0,$$
then $D$ is nef.
\end{corollary}

\begin{proof}
First using Mumford's identity on $\M_g$ 
$$-\delta_{0}=-12\lambda+ \kappa_1 + \sum_{1 \le i \le \lfloor \frac{g}{2} \rfloor}\delta_i,$$
write 
$$D= (a-12b_{0}) \lambda + b_{0} \kappa_1 + \sum_{1 \le i \le \lfloor \frac{g}{2} \rfloor}(b_{0}-b_i)\delta_i.$$
Therefore by Lemma  \ref{pullback}, one has

\begin{eqnarray*}
f^*D = b_{0} \kappa_1 + (b_1 - b_{0})\sum_{1 \le i \le g} \psi_i 
+ \sum_{2 \le i \le \lfloor \frac{g}{2} \rfloor}(b_{0}-b_{i}) B_i.
\end{eqnarray*}

Substituting the relation $\kappa_1 = K_{\M_{0,g}}+\sum_{}B_i$, and then big averaging the $\psi_i$, one has
\begin{eqnarray*}
f^*D = b_{0} K_{\M_{0,g}} + (b_1 - b_{0})\sum_{1 \le i \le g} \psi_i 
+ \sum_{2 \le i \le \lfloor \frac{g}{2} \rfloor}(2b_{0}-b_{i}) B_i.\\
=b_0 K_{\M_{0,g}} 
+ \sum_{2 \le i \le \lfloor \frac{g}{2} \rfloor}(\frac{i(g-i)}{(g-1)}(b_1-b_0) + 2b_{0}-b_{i}) B_i.\\
\end{eqnarray*}
It is enough to show that under the given hypothesis, the coefficients of the $B_i$ above
are nonnegative and  $\le b_0$, so that by the Ray Theorem,
$f^*D$ and hence $D$ is nef.  That is,  for $2 \le i \le \lfloor \frac{g}{2} \rfloor$, 
$$-b_0 (g-1) \le i(g-i)b_1+(i^2-ig+g-1)b_0-(g-1)b_i=i(g-i)(b_1-b_0)+(g-1)(b_0-b_i) \le 0.$$
By hypothesis, this is true. 
\end{proof}

\begin{corollary}\label{level0criteria}\label{level0filter} Let $D= a\lambda-\sum_{0\le i \le \lfloor \frac{g}{2} \rfloor}b_i \delta_i$ be
an $F-$divisor on $\M_{g}$.  If there exists a constant $c \ge 0$ such that 
$$\frac{2g-2-i(g-i)}{g-1}+\frac{b_1i(g-i)-b_i(g-1)}{(g-1)c} \le c,$$
for all $i \in \{2 \ldots \lfloor \frac{g}{2} \rfloor\}$, then $D$ is nef.
\end{corollary}

\begin{proof}
By The Bridge Theorem, $D$ is nef as long as $f^*D$ is nef.  To show the assumptions in the theorem guarantee that
$f^*D$ is nef use the proof of Theorem \ref{step0} with the Ray Theorem.
\end{proof}

\noindent
It may be interesting to note that Corollary \ref{level0criteria} can't seem to be improved using Mumford's criteria.

\bigskip

\noindent
A divisor that doesn't meet the conditions above can of course 
still be nef.  For example, in Corollary \ref{level0criteria}, if whatever constant $c$ is tried, there is a boundary class in the support of $D$ with a coefficient
larger than $c$, then more has to be done to show $D$ is nef.   In particular, one can still prove $D$ is nef by showing the divisor is nef when restricted to the boundary component whose class
has coefficient bigger than $c$.  By assuming more about the divisor, say that every boundary restriction  has to be nef, 
one obtains the criteria in the next two results.  The first, Corollary \ref{$b_1$} comes from using Theorem \ref{sometimesrestrict} with $c=0$. 	The remaining criteria of the section are Corollaries of this fact.  Each provides an easy to check condition which guarantees a divisor on $\M_g$ is nef.

\begin{corollary}\label{$b_1$}An F divisor $D =a \lambda  -\sum_{i=0}^{\lfloor \frac{g}{2} \rfloor}b_i \delta_i$ 
 on $\M_{g}$ is nef provided that  $b_i \le b_1$ for all $i \ge 2$.
 \end{corollary}

	\begin{proof}
	Let $D$ be as described in the hypothesis of the theorem.  It will be shown that any boundary restriction of an F divisor $f^*D$
	is equivalent to an effective sum of boundary classes.  For simplicity of notation, put $D=f^*D$.

	Let $v=v_{a,z}:\M_{0,a+z}\longrightarrow \M_{0,g}$ be a boundary restriction where we attach   
	an $n_i+1 \ge 3$ pointed curve to each point $p_i \in A$, where $|A|=a$ and
	do nothing to the $z$ points $q_i \in Z$.  Then as we have seen in Lemma $1$,
	$$v^*D=b_1\sum_{i\in Z=\{i|n_i=2\}}\psi_{i}+\sum_{i\in A=\{i|n_i>2\}}b_{n_i}\psi_{i}-\sum_{\stackrel{B \subset A, y\le |Z|}{2\le y+|B|\le \lfloor \frac{g}{2}\rfloor}}b_{y+\sum_{i\in B}n_i}\Delta^{Z,y}_{A,B},$$ where $\Delta^{Z,y}_{A,B}=\sum_{Y \subset  Z, |Y|=y} \delta_{Y \cup B}$.

	The proof is divided into 3 cases: $z \ge 4$, $z=3$ and $z=2$.  First suppose that $z \ge 4$.  Let $Z = \{1,\ldots,z\}$.  By averaging the $\psi_i$, the divisor $$D_z=\sum_{i \in Z}\psi_{i}-\sum_{\stackrel{S \subset Z}{2 \le s=|S| \le \lfloor \frac{z}{2} \rfloor}}\delta_{S} = \sum_{\stackrel{S \subset Z}{2 \le s=|S| \le \lfloor \frac{z}{2} \rfloor}}(\frac{(s-1)z-s^2+1}{z-1})\delta_{S}$$ in $\M_{0,z}$.  Each coefficient is positive as long as $z \ge 4$.  To see this, put $f(s)=(s-1)z-s^2+1$.  Then $f'(s)=z-2s \ge 0$ since $s \le \frac{z}{2}$.  So the function $f(s)$ is increasing in the range we are interested in.  Now as $f(2) \ge 1$, $f$ is always positive. In particular, $\pi_a^*(D_z) $
	is an effective sum of boundary classes in $\M_{0,a+z}$.
	Now let $\pi_a:\M_{0,a+z}\longrightarrow \M_{0,z}$ be the morphism which drops the attaching points $p_i \in A$.  Then 
	$$v^*D -b_1 \pi_a^*(D_z)   =  
	\sum_{i \in A}b_{n_i} \psi_i-
	\sum_{I \subset A}b_{\sum_{i\in I}n_i} \Delta^{Z,0}_{A,I}
	+\sum_{\stackrel{y>0, I \subset A}{0 \le |I| \le a}}(b_1-b_{y+\sum_{i\in I}n_i}) \Delta^{Z,y}_{A,I}.$$
	For $y >0$, the coefficients of the classes $\Delta^{Z,y}_{A,I}$ are nonnegative since by hypothesis, $b_1 \ge b_i$ for all $i$.  Fix two elements $p,q \in Z$.  Then for $i \in A$, $\psi_i = \sum_{I \subset \{p,q\}^c}\delta_{I\cup i}$ and so $$\sum_{i \in A}b_{n_i}\psi_i-\sum_{I \subset A}b_{\sum_{i\in I}n_i}\Delta^{Z,0}_{A,I}=\sum_{I \subset A} (\sum_{i\in I}b_{n_i} - b_{\sum_{i\in I}n_i})\Delta^{Z,0}_{A,I} + E,$$ where E is an effective sum of boundary classes.  That the coefficients $(\sum_{i\in I}b_{n_i }- b_{\sum_{i\in I}n_i}) \ge 0$ is a consequence of the assumption that $D$ is an F divisor and so it's coefficients satisfy property $5$ of Thm. $1$.

	Now suppose that $z=3$.  By replacing the $\psi_i$ for $i \in Z$ by their averages and by using the same partial average for the $\psi_i$ for $i \in A$ as was done in the previous case, we get that:
	\begin{multline}
	$$v^*D=
	\sum_{\stackrel{B \subset A,|B|=b}{2\le b\le a}}(\frac{b_1 3b(b-1)}{(a+1)(a+2)}+(\sum_{i \in B}b_{n_i}-b_{\sum_{i\in B}n_i}))\Delta^{Z,0}_{A,B} 
	+\\
	\sum_{\stackrel{B \subset A, |B|=b}{1\le b\le a-1}}
	(\frac{b_1 2(2+b)(1+b)}{(a+2)(a+1)}+(b_1-b_{1+\sum_{i\in B}n_i}))\Delta^{Z,y}_{A,B}.\\
	\end{multline}
	These coefficients are nonnegative since by hypothesis $b_1 \ge b_i \ge 0$ for all $i$ and since by assumption $D$ is an F divisor and so by Proposition \ref{F-inequalities}, $(\sum_{i\in I}b_{n_i }- b_{\sum_{i\in I}n_i}) \ge 0$.

	Consider the case $z=2$.  For $p \in Z$, we form a partial average of $\psi_p$ by taking $q=Z \setminus p$ and fixing any $i \in A$ so that $\psi_p = \sum_{I \subset A\setminus \{i\}}\delta_{I\cup p}.$
	There are $a$ ways of fixing such a point $i \in A$.  So $a \psi_p = \sum_{I \subset A,|I|=i} (a-i)\delta_{I\cup p}$ and so for $Z = \{p,q\}$:
	$$(\psi_p + \psi_q)=\sum_{\stackrel{I \subset A,|I|=i}{1\le i \le a-1}}(\frac{a-i}{a}+\frac{a-(a-i)}{a})\Delta^{Z,1}_{A,I}=\sum_{\stackrel{I \subset A}{1\le |I| \le a-1}}\Delta^{Z,1}_{A,I}.$$

	Once again, by replacing the  $\psi_i$ for $i \in A$ as was done in the two previous cases, we get that:
	$$v^*D=
	\sum_{\stackrel{B \subset A}{2\le |B|\le a}}(\sum_{i \in B}b_{n_i}-b_{\sum_{i\in B}n_i})\Delta^{Z,0}_{A,B} 
	+
	\sum_{\stackrel{B \subset A}{1\le |B|\le a-1}}
	(b_1-b_{1+\sum_{i\in B}n_i})\Delta^{Z,y}_{A,B}.$$

	These coefficients are nonnegative by assumption.  Therefore, any F divisor $D=b_1\sum_{i=1}^{g}\psi_i-\sum_{i=2}^{\lfloor \frac{g}{2} \rfloor}b_iB_i$ in $\M_{0,g}$ such that $b_i  \le b_1$ for all $i$ is nef.
\end{proof}

\begin{corollary} \label{$b_m$}An F divisor $D =a \lambda  -\sum_{i=0}^{\lfloor \frac{g}{2} \rfloor}b_i \delta_i$ 
 on $\M_{g}$ is nef provided that  $$2 \ min \{b_i\ |\ i \ge 1\} \ge max \{b_i\ |\ i \ge 1\}.$$ 
\end{corollary}
			
\begin{proof}
Let $D$ be as described in the hypothesis of the theorem.  It will be shown that any boundary restriction of an F divisor $f^*D$
	is equivalent to $cK +E$, where $E$ is an effective sum of boundary classes for some $c \ge 0$.  For simplicity of notation, put $D=f^*D$.

						Let $v=v_{a,z}:\M_{0,a+z}\longrightarrow \M_{0,g}$ be a boundary restriction where we attach   
							an $n_i+1 \ge 3$ pointed curve to each point $p_i \in A$, where $|A|=a$ and
								do nothing to the $z$ points $q_i \in Z$.  Then as we have seen in Lemma $1$, 
								for $Z=\{i|n_i=2\}$ and $A=\{i|n_i>2\}$:
								
\begin{multline*}
v^*D=b_1\sum_{i\in Z }\psi_{i}+\sum_{i\in A}b_{n_i}\psi_{i}-\sum_{\stackrel{B \subset A, y\le |Z|}{2\le y+|B|\le \lfloor \frac{g}{2}\rfloor}}b_{y+\sum_{i\in B}n_i}\Delta^{Z,y}_{A,B} \\
= (b_1-c)\sum_{i\in Z}\psi_{i}+\sum_{i\in A}(b_{n_i}-c)\psi_{i}
+c\sum_{i\in A \cup Z}\psi_i-\sum_{\stackrel{B \subset A, y\le |Z|}{2\le y+|B|\le \lfloor \frac{g}{2}\rfloor}}b_{y+\sum_{i\in B}n_i}\Delta^{Z,y}_{A,B} \\
=(b_1-c)\sum_{i\in Z}\psi_{i}+\sum_{i\in A}(b_{n_i}-c)\psi_{i}
+cK_{\M_{0,a+z}} + \sum_{\stackrel{B \subset A, y\le |Z|}{2\le y+|B|\le \lfloor \frac{g}{2}\rfloor}}(2c-b_{y+\sum_{i\in B}n_i})\Delta^{Z,y}_{A,B},
\end{multline*}
where $\Delta^{Z,y}_{A,B}=\sum_{Y \subset  Z, |Y|=y} \delta_{Y \cup B}$.   Recall, as is explained in Section $2$, 
each class $\psi_i$ is equivalent to an effective sum of boundary classes.  So as long as $c \le b_i \le 2c$ for all $i$, then 
$v^*D = cK_{\M_{0,a+z}}+E$ as required.  Just take $c \in [max\{b_i \ | \ i \ge 1\}/2,min\{b_i \ | \ i \ge 1\}]$,
which, by hypothesis, is a nonempty interval.
\end{proof}

\begin{corollary}\label{induct}
Let $D= a \lambda-\sum_{i=0}^{\lfloor \frac{g}{2} \rfloor}b_i \delta_i$ be an $F-$divisor on $\M_g$.  If $g$ is odd and $b_j=0$ 
or if $g$ is even and $b_j=0$ for $j < \frac{g}{2}$, then $D$  is nef.  
\end{corollary}

To prove Corollary \ref{induct}, the following result will be used.

\begin{lemma}\label{technical}
If  $D = a \lambda -\sum_{i=0}^{\lfloor \frac{g}{2} \rfloor}b_i \delta_i$ is an F-divisor in $\M_g$ such that  $b_i=0$ , then 
\begin{enumerate}
\item$b_j=b_k$  for all $j,k$ such that $j+k=i$, and 
\item $b_{j}=b_{i+j}$ for all $j \ge1$ such that $i+j \le g-1.$
\end{enumerate}
\end{lemma}

\begin{proof}Since $D$ is an F divisor, $b_{g-(j+k)}+b_{j}-b_k \ge 0$ and  $b_{g-(j+k)}+b_k-b_j \ge 0$.  But $b_{g-(j+k)}=b_i=0$ and so $b_j=b_k$.  The second assertion follows from the fourth type of inequality $b_{i}+b_{j}\ge b_{i+j}$, which since $b_i=0$, gives that  $b_{j}\ge b_{i+j}$.  By substituting $b_{j}=b_{g-j}$ and $b_{i+j}=b_{g-(i+j)}$,  one has
$b_{i+j}\ge b_{j}$. 
\end{proof}

\begin{proof}(of Corollary \ref{induct}) Let  $D = a \lambda -\sum_{i=0}^{\lfloor \frac{g}{2} \rfloor}b_i \delta_i$ in $\M_g$ be  an F-divisor such that $b_j=0$ for some $j$.  The result will be proved by induction on $j$.
Of course if $b_1=0$, the divisor is trivial and there is nothing
to prove.  If $b_2=0$, then by Lemma \ref{technical}, $b_2=b_{2x}=0$ for all $x$ such that $2x \le g-1$, and $b_1=b_{1+2x}$
for all $x$ such that $1+2x \le g-1$.  Therefore, $b_i \le b_1$ for all $i$ so that by Corollary \ref{$b_1$}, $D$ is nef.

Suppose $b_k = 0$ for some $3 \le k < \lfloor \frac{g}{2} \lfloor$ and that  the statement is true when $b_i=0$ for all $i <k$.   Consider $m$ so that $mk \le g-1$ but that $(m+1)k > g-1$. By Lemma \ref{technical}, $0=b_k=b_{mk}=b_{g-mk}$.  Then $g-mk < k$, and so $b_{g-mk}=0$ means that by induction, the statement is true.   

Now suppose that $g=2n-1$ is odd and $b_{ \lfloor \frac{g}{2} \rfloor}=b_n=0$.  Then
by Lemma \ref{technical}, $b_{n}=b_{2n}=b_{1}=0$.  Hence $b_i =0$ for all $i\ge1$ and $D$ satisfies Corollary \ref{$b_1$}.
\end{proof}

\section{Using the Nef Wizard to show that the criteria prove Conjecture $1$ for low values of $g$}

One can show by a computer check that all the $F-$divisors in $\M_g$, 
for at least  $g \le 24$ are nef. 

\begin{theorem}\label{22}
The $F-$conjecture is true on $\M_{0,g}/S_g$ for $g \le 24$.
\end{theorem}

\begin{corollary} \label{F22}
The $F-$conjecture is true on $\M_{g}$ for $g \le 24$.
\end{corollary}

\begin{proof} (of Corollary \ref{F22})
Apply \cite{GKM} Theorem $.7$.
\end{proof}

The procedure for doing so is explained in this section.
The starting point is that by \cite{GKM}, the conjecture on $\M_g$ is 
equivalent to the conjecture on $\M_{0,g}/S_g$.  In particular, if one can prove that the extremal $F-$divisors
on $\M_{0,g}/S_g$ are nef, then the $F-$conjecture is true on $\M_g$.  The computer program Nef Wizard 
generates the extremal $F-$divisors on $\M_{0,g}/S_g$ in terms of the sums 
of boundary classes $B_i$.  Nef Wizard finds $F-$divisors on $\M_g$ that pullback to the extremal divisors via $f$ so that
the criteria may be applied.

To prove Theorem \ref{22},  the following result will be used.

\begin{lemma}\label{E}
Let $E=\sum_{2 \le i \le \lfloor \frac{g}{2} \rfloor}e_i B_i$ be a divisor on $\M_{0,g}/S_g$ and
consider $$D_E=a\lambda - b_0 \delta_0 - b_1 \delta_1 
- \sum_{2 \le i \le \lfloor \frac{g}{2} \rfloor}(\frac{i(g-i)}{(g-1)}b_1-e_i)
		\delta_{i},$$
		where 
		\begin{enumerate}
		\item $b_1=max\{0,\frac{(g-1)}{i(g-i)} e_i,\frac{(g-1)}{2ij}(e_i+e_j-e_{i+j}) \ | \ 1 \le i,j; I+j \le g-1  \}$, \\
		\item $b_0=\frac{1}{2}max\{b_i \ | \ i \ge 1 \}$, and \\
		\item $a=12b_0 - b_1$.
		\end{enumerate}
		Then $f^*D_E=E$ and If $E$ is an $F-$divisor then so is $D_E$.  
		\end{lemma}

		\begin{proof}To see that $f^*D_E=E$ use Lemma's $1$ and $5$: 
		\begin{multline*}
		f^*D_E= b_1 \sum_{1 \le i \le g}\psi_i - \sum_{2 \le i \le \lfloor \frac{g}{2} \rfloor}(\frac{i(g-i)}{(g-1)}b_1-e_i)B_i\\
			=\sum_{2 \le i \le \lfloor \frac{g}{2} \rfloor}(\frac{i(g-i)}{(g-1)}b_1-(\frac{i(g-i)}{(g-1)}b_1-e_i))B_i.
			\end{multline*}
			Now suppose that $E$ is an $F-$divisor. To show that $D_E$ is also an $F-$divisor one just checks that it satisfies the five
			inequalities of Theorem \ref{F-inequalities}.  The first four are true by definition of $D_E$.  For example,
			to see that $b_i+b_j - b_{i+j} \ge 0$:
			\begin{multline*}
			(\frac{i(g-i)}{(g-1)}b_1-e_i)+(\frac{j(g-j)}{(g-1)}b_1-e_j)-(\frac{(i+j)(g-(i+j))}{(g-1)}b_1-e_{i+j})\\
			=\frac{2ijb_1}{(g-1)}-(e_i+e_j-e_{i+j}),
			\end{multline*}
			which is nonnegative as long as $$b_1 \ge \frac{(g-1)}{2ij}(e_i+e_j-e_{i+j}).$$
			The fifth inequality holds because $f^*D_E=E$.
			\end{proof}

\begin{proof} (of Theorem \ref{22})
By using a computer program such as LRS \cite{lrs}, one can generate a list of extremal divisors $E$ for the $F-$cone of $\M_{0,g}/S_g$.  This computation is convenient to perform by considering divisors expressed in the basis for Pic$\M_{0,g}/S_g$ given by
 $\{B_i\}_{2 \le i \le \lfloor \frac{g}{2} \rfloor}$.   To change these extremal divisors into the form necessary to apply the theorems, one can solve for $D_E$ as in Lemma \ref{E} and then pull back.  Finally, to check the divisors are were all nef, we ran them 
 through the
 program Nef Wizzard. 
 \end{proof}

\section{Relevance of Conjecture}
If the $F-$conjecture is true it means that  the extremal rays of the cone of curves $\overline{NE}(\M_{g})$ are spanned by
the $F-$curves.  This would be very good information to have since as was illustrated in the Introduction, $\overline{NE}(\M_{g})$
reveals  information about the birational
 geometry of $\M_{g}$.  Moreover, it would mean that $\overline{NE}(\M_{g})$  is an interesting example of a cone of curves.  
 To explain why, I'll say a little bit about he minimal model program (MMP). 
 
 The MMP generalizes the birational classification of smooth 
	 surfaces using certain kinds of projective
	 morphisms called contractions.  Contractions are morphisms $f:X \longrightarrow Y$ between projective varieties such that
	 $f_*(\mathcal{O}_X)=\mathcal{O}_Y$; they are determined by the faces  of the cone of curves. 
	 Unlike the situation for surfaces, for higher dimensional projective varieties, contractions are not so resolutely understood 
	 nor is their existence guaranteed.  

	 In order  to classify $X$ using contractions $X \longrightarrow Y$ one studies the image variety $Y$ and the fibers of the contraction morphism. 
	 	There are a couple
					 of possibilities depending on whether or not the image $Y$ has the same dimension as $X$. 
					 			 If $\mbox{dim}X > \mbox{dim}Y$, this is a fibral type contraction.  As was mentioned above, by \cite{thesis},
								 			 there are no fibral 
											 			 type contractions of $\M_{g}$.   The other possibility is that $\mbox{dim}X = \mbox{dim}Y$.  For $\M_{g}$ and other 
														 			 higher dimensional varieties $X$, two things can happen.  The first is that the morphism $X \longrightarrow Y$
																	 			 is a so called divisorial 
																				 			 contraction -- this is the analog of the surface case wherein $X$ is the blowup of $Y$.  
																							 			 By Proposition $6.4$ in \cite{GKM},  for $g \ge 5$ the only divisorial contraction of $\M_g$ is a blowdown of elliptic tails.      When $g=3$ there is another divisorial contraction (\cite{Rulla}) and the problem is open when $g=4$.  The remaining
																										 			 kind of contraction  doesn't have an analog in the classification of surfaces.  It is called a small contraction 
																													 			 and it is																																 essentially the case where the image variety $Y$ has bad singularities and so one has to surgically repair
																																 			 it (i.e.  do flips or flops) in order to proceed with the program.

																																			 As stated in the introduction, since there are a  finite number $F-$curves to begin with, if the $F-$Conjecture is true, the cone of curves is polyhedral, like the cone of curves for a Fano variety.  This 
																																			 is counter-intuitive since for $g \ge 23$, the Kodaira Dimension
																																			 of $\M_g$ is positive (in fact, for $g \ge 24$, the moduli space is of general type).  
																																																																																																											Finally, when one considers $\M_g$ to be defined over a field of positive characteristic, then every extremal face of $\overline{NE}(\M_{g})$ gets contracted.  This is also surprising since contractions of a variety $X$ are only guaranteed for $K_X-$negative extremal rays, and only one of the F-curves is $K_{\M_g}$ negative.
																																																																																																																																			
																																																																																																																																						In any case, for low genus when  the nef cones and the $F-$cones of $\M_g$ are the same,
																																																																																																																																									one has a series of explicit examples	of cones of curves that 
																																																																																																																																												have finitely many extremal rays,  each spanned by a smooth, irreducible and 
																																																																																																																																															rational curve.  Moreover, when the characteristic of the field is positive,
																																																																																																																																																		every face of the cones get contracted,  none	of the contractions is fibral, and in fact all
																																																																																																																																																					but one  are small contractions.   Hence one has a collection of rich examples which
																																																																																																																																																											deepen our understanding of the birational geometry of the spaces $\M_{g}$.  Moreover,
																																																																																																																																																																	though admittedly not the simplest of examples, these cones broaden our understanding of 
																																																																																																																																																																										cones of curves in general.

\nocite{keelthesis}
\nocite{MoKo}
\nocite{Har}
\nocite{lrs}
\nocite{CR}
\nocite{149}
\nocite{danny}

\bibliographystyle{alpha}
\bibliography{citations}

\vskip 5pt\noindent DEPARTMENT OF MATHEMATICS, UNIVERSITY OF MICHIGAN, ANN ARBOR, MI 48109-1109\vskip 3pt\mbox{ } E-mail:\ {\tt agibney@umich.edu}

\end{document}